\documentclass[10pt, twoside, reqno]{amsart}

\usepackage{amssymb}
\usepackage{amsmath}
\usepackage{amsthm}
\usepackage{amsopn}
\usepackage{mathrsfs}

\theoremstyle{plain}
\newtheorem{theorem}{Theorem}
\newtheorem{lemma}{Lemma}
\newtheorem{prop}{Proposition}

\theoremstyle{definition}

\newtheorem{remark}{Remark}
\newtheorem{corollary}{Corollary}

\def\arg{\mathop{\rm arg}}
\def\argmax{\mathop{\rm arg\max}}

\def\esp{{\bf E}}
\def\Pb{{\bf P}}

\def\RR{\mathbb R}

\def\NN{\mathbb N}
\def\Ln{\textit{l}_n}
\def\pln{\textsl{p}\textit{l}_n}
\def \1{{\rm 1}\mskip -4,5mu{\rm l} } 

\def\eps{\varepsilon}

\def\theta{\vartheta}

\def\ds{\displaystyle}
\def \1{{\rm 1}\mskip -4,5mu{\rm l}}
\def\hkapjsq{{h_{\kappa_j}^*\!\!}^2}

\begin{document}

\parskip5pt
\parindent0pt

\title[Semiparametric Shift Estimation]{Stein Shrinkage and Second-Order Efficiency for semiparametric estimation of the shift}
\date\today


\author[A. S. Dalalyan]{Arnak S.\ Dalalyan}

\address{\vbox{Universit\'e Paris 6\\
Laboratoire de Probabilit\'es\\
Bo\^\i te courrier 188\\
75252 Paris Cedex 05, France\\
\texttt{dalalyan@ccr.jussieu.fr}
}}

\begin{abstract} The problem of estimating the shift (or, equivalently, the center of symmetry)
of an unknown symmetric and periodic function $f$ observed in Gaussian white noise is
considered. Using the blockwise Stein method, a
penalized profile likelihood with a data-driven penalization is introduced so that
the estimator of the center of symmetry is defined as the
maximizer of the penalized profile likelihood. This estimator
has the advantage of being independent of the functional class to which the signal $f$ si assumed to belong
and, furthermore, is shown to be semiparametrically adaptive and efficient.

Moreover, the second-order term of the risk expansion of the proposed estimator is proved to behave at least
as well as the second-order term of the risk of the best possible estimator
using monotone smoothing filter. Under mild assumptions, this
estimator is shown to be second-order minimax sharp adaptive over
the whole scale of Sobolev balls with smoothness $\beta>1$. Thus,
these results extend those of \cite{DGT}, where second-order asymptotic minimaxity is proved for an
estimator depending on the functional class containing $f$ and $\beta\geq 2$ is required.
\end{abstract}



\maketitle

\section{Introduction}\label{sec1}
\subsection{Motivation}
The increasing interest to the semiparametric models in
recent years is mainly explained by the compromise they offer between the
relative simplicity of parametric inference and the flexibility of nonparametric
modeling. In many semiparametric models, though the ignorance of an
infinite dimensional nuisance parameter, the finite-dimensional parameter of
interest can be estimated as well as if the nuisance parameter were known.
In most situations, there are many estimators having this feature. Moreover,
most of them depend on the functional class to which the unknown nuisance
parameter is assumed to belong.

The aim of the present paper is to propose a second-order efficient and
entirely data-dependent estimator in the problem of shift estimation
when the observed shifted signal is corrupted by the Gaussian white noise.
This problem has been studied in \cite{DGT},
and the goal of the present paper is to complete the inference developed
there. In particular, the estimator proposed in \cite{DGT} and proved to be
second-order minimax when the signal belongs to a Sobolev ball, depends on the
parameters of the Sobolev ball in question. In the present work, we propose an
estimator of the shift parameter that is independent of the parameters of the Sobolev ball
and is second-order minimax simultaneously for a large scale of Sobolev balls.

To further motivate the study of the second-order efficiency in semiparametric
inference, let us briefly recall a popular general method of construction of efficient
estimators in a semiparametric model: the profile likelihood maximization.
The asymptotic properties of this method are studied in \cite{SW} and, in a more general fashion,
in \cite{MuVaa}.
Profile likelihood techniques are effectively applied in a number of contexts
such as laser vibrometry signals \cite{LLL}, varying
coefficient partially linear models \cite{FanHuang}, and so forth.

If the parameter of interest is partitioned as $(\vartheta; f)$, with $\vartheta$ being a
low-dimensional parameter of interest and $f$ a higher dimensional nuisance parameter, and $\Ln(\vartheta; f)$
is the log-likelihood of the model, then the profile likelihood for $\vartheta$ is defined as
$\pln(\vartheta) =\sup_{f\in\mathcal F} \Ln(\vartheta; f)$ and the Profile Likelihood
Estimator (PLE) is $\vartheta_{PLE} = \argmax_\vartheta \pln(\vartheta)$. Thus, the nuisance parameter $f$
is eliminated by taking the sup over all possible values of $f$ in some
a priori chosen class $\mathcal F$.

A natural question arises: what is the best way of choosing the class F and
what is the impact of this choice on the accuracy of the PLE? The theory
fails to answer this question as long as only the first-order term of the risk is
considered. It seems that the most appealing way to study the dependence
on $\mathcal F$ of the accuracy of the PLE is to consider the second-order term of the
quadratic risk. This approach is developed in \cite{GH1,GH2} for partial linear models,
in \cite{DGT} for a nonlinear model with a shift parameter, and in \cite{Cast} for a
model with scaling parameter. In these papers, the second-order asymptotic minimaxity of the
proposed estimators is proved and optimal constants are obtained.

Another important point is that the study of the second-order asymptotics of the risk allows one to refine the choice
of the tuning parameters, if there are, of the estimating procedure. This argument is used
in \cite{ht,Kang,MamPark} in order to propose a rate-optimal choice of a tuning parameter.

The results of the above mentioned papers grant an increasing importance to the
second-order terms in that they show that, in a semiparametric estimation problem,
the second-order term is not dramatically smaller than the first-order term, especially
when the nuisance parameter is not very smooth (or not very sparse). Thus,
the investigation of second-order efficient estimators is not only a
challenging theoretical problem, but  is also of practical interest.

\subsection{Methodology} The main goal of the present paper is to define an estimator
of $\vartheta$ which adapts automatically to the smoothness of $f$
and is simultaneously second-order efficient over a large variety of Sobolev
balls. For linear models, such a procedure has firstly been proposed by Golubev and
H\"ardle in \cite{GH2}. They use the well known idea of unbiased risk minimization
in order to determine the data driven filter. However, their procedure
is provably second-order efficient only when the (data driven) filter and the contrast function are based
on independent sub-samples. This sample-splitting technique is
frequently used in theory (see e.g.\ \cite{B}, \cite[p.\ 19]{pf}), but it
is rather unattractive from the practical point of view.

In this paper, we focus on the shift estimation of a periodic symmetric
signal and use the celebrated Stein shrinkage (see \cite{JS,Ste}) for defining the
data-driven version of the penalized profile likelihood. While there is a huge number of papers
concerning the Stein's shrinkage in nonparametric estimation
(\cite{BLZ, Cai, CaiLow, CavTsy2, CavTsy1, DJ,RIG}, see also \cite{TSY}
for a comprehensive discussion), this is to our knowledge the first paper where
Stein's shrinkage is applied in the context of semiparametric estimation.

The construction we use is closer to
the one of \cite{CavTsy1} in that a penalized version of the
Stein estimator with weakly geometrically increasing blocks is considered.
It seems that, unlike in nonparametric inference, in semiparametric inference
it is crucial to use the penalized shrinkage in order to get the second-order efficiency.
An important advantage of using the Stein Shrinkage is that, unlike the procedure proposed in \cite{GH2},
our procedure is not based on the sample-splitting technique and, nevertheless, enjoys the property of
second-order asymptotic minimaxity.

Note also that our results hold when the function $f$ has arbitrary Sobolev smoothness $\beta>1$. However,
for $\beta<2$, the penalized profile likelihood (PPL) is not necessarily concave over the whole
parameter space, therefore we use a two-step procedure. We first construct an initial $\sqrt n$-consistant
estimator $\bar\theta$ of $\theta$, and then maximize the PPL over a shrinking
neighborhood of $\bar\theta$, where we prove that the PPL is concave with probability close to one.

\subsection{Structure of the paper} Section 2 describes the model and introduces
the Penalized Maximum Likelihood Estimator (PMLE) based on a filtering sequence $h$.
In Section 3, the local concavity of the PPL is proved and the PMLE based on a data dependent choice of
$h$ is proposed. Oracle inequalities for adaptive PMLE and its second-order
efficiency over Sobolev balls are stated and proved in Section 4. Finally, Section
5 contains the definition of a preliminary estimator and the technical
details of the proofs.

\section{A simple  semiparametric model}\label{S2}

Consider the ``signal in Gaussian white noise model'', that is the
observations $(x^\eps(t),\,t\in[-1/2,1/2])$ with
\begin{equation}\label{x^eps}
dx^\eps(t)=f_\theta(t)\,dt+\eps dW(t),\qquad t\in[-1/2,1/2],
\end{equation}
are available, where $W(t)$ is a Brownian motion. Assume that the
signal has the form $f_\theta(t)=f(t-\theta)$, where $f:\RR\to\RR$ is
a symmetric periodic function having $1$ as smallest period. More precisely, we assume that
the function $f$ belongs to the set $\mathcal F_0=\cup_{\rho>0}\mathcal F_0(\rho)$ with
$$
\mathcal F_0(\rho)=\Big\{f\in L^2_{loc}\,:\,
f(x)=f(-x)=f(x+1),\ \forall\,x\in\RR;\ |f_1|\geq \rho\Big\},
$$
where we denote by $L^2_{loc}$ is the set of all locally squared integrable functions and
by $f_1=\int_{-1/2}^{1/2}f(t)\,\cos(2\pi t)\,dt$.

The goal is to estimate the parameter $\theta\in\Theta\subset]-T,T]$
with $T<1/4$. As explained in \cite{DGT}, the assumption $T<1/4$ is necessary for the identifiability
of the parameter $\theta$. In this context, the unknown function $f$ is considered
as an infinite dimensional nuisance parameter.

The first-order asymptotic properties of
estimators in closely related models have been studied in \cite{G, hm, Sch1,Sch2}. Note also that the model
we have just defined is an idealized version of the symmetric location model
\cite{Sto, MamPark} and the shifted curves model \cite{GLM}.

Let us introduce the ``sequence model'' equivalent to (\ref{x^eps}). To do this,
for any integer $k\geq 1$, we denote
\begin{align}\label{xx*}
\begin{matrix}
x_k&=\sqrt 2\,\hbox{$\ds\int_{-1/2}^{1/2} \cos(2\pi k t)\,dx^\eps(t)$},\\
x_k^*&\ds=\sqrt 2\,\hbox{$\ds\int_{-1/2}^{1/2} \sin(2\pi k t)\,dx^\eps(t)$}.
\end{matrix}
\end{align}
Clearly,
\begin{eqnarray}
\begin{matrix}\label{x_k}
x_k&=f_k\cos(2\pi k\theta)+\eps \xi_k,\\
x_k^*&=f_k\sin(2\pi k\theta)+\eps \xi_k^*,
\end{matrix}
\end{eqnarray}
where $f_k=\sqrt 2\int_{-1/2}^{1/2} \cos(2\pi k t)\,f(t)\,dt$ and
$(\xi_k,\xi_k^*,k=1,2,\ldots)$ are independent standard Gaussian
random variables. Note that, the laws of likelihood processes (indexed by
$(\theta,f)\in\RR\times \mathcal F_0$) of the models (\ref{x^eps})
and (\ref{x_k}) coincide implying thus the equivalence of these
models. It is useful to mention here that the Fisher information in the problem of estimating $\theta$ with fixed $f$
is
$$I^\eps(f)=\eps^{-2}\int_{-1/2}^{1/2}f'(x)^2\,dx=\eps^{-2}\sum_{k\in\NN} (2\pi k)^2 f_k^2.$$
In this paper, we estimate the parameter $\theta$ by a version of
the method of profile likelihood maximization (cf.\ \cite[p.\ 106]{vdv3}), which is also known as
penalized maximum likelihood estimator \cite{DGT}. We recall briefly its definition. Write
$\Pb_{\theta,f}$ (resp.\ $\esp_{\theta,f}$) for probability measure
(resp.\ expectation) induced by $x^\eps$ on the canonical space
$C([-1/2,1/2])$ equipped with the Wiener measure. As no confusion is
possible, we use the same notation in the ``sequence model'' given
by (\ref{x_k}). The Radon density $L_\eps(\tau,f,\cdot)$ of
$\Pb_{\tau,f}$ with respect to $\Pb_{0,0}$ is given by
\begin{align*}
L_\eps(\tau,f,x^\eps)&=\exp\bigg(\eps^{-2}\int_{-1/2}^{1/2}f(t-\tau)\,dx^\eps(t)-\frac{\eps^{-2}}{2}
\int_{-1/2}^{1/2}f^2(t)\,dt\bigg)\\
&=\exp\bigg(\frac{\sqrt2}{\eps^{2}}\sum_{k=1}^\infty
f_k\int_{-1/2}^{1/2}\cos[2\pi
k(t-\tau)]\,dx^\eps(t)-\frac{1}{2\eps^{2}} \sum_{k=1}^\infty
f_k^2\bigg).
\end{align*}
Easy algebra yields
\begin{align}\label{ML}
\max_{(f_k)_k\in \RR^\NN}
L_\eps(\tau,f,x^\eps)&=\exp\bigg\{{\eps^{-2}}\sum_{k=1}^\infty
\bigg(\int_{-1/2}^{1/2}\cos[2\pi
k(t-\tau)]\,dx^\eps(t)\bigg)^2\bigg\}.
\end{align}
For any $\tau\in\RR$, this expression is equal to infinity for
almost all paths $x^\eps$. Thus, it is necessary to impose some restrictions on the set over
which the maximization is done.  This is equivalent to considering a penalized profile
likelihood. In this paper, we focus on the penalization corresponding to
ellipsoids in $\ell^2$. More details on this method can be found in
\cite{DGT}, we here content ourselves with giving the final definition.

We call {\it filtering sequence} or \textit{filter} any
$h=(h_k)_{k\in\NN}\in[0,1]^{\NN}$ such that only a finite number of
$h_k$ are non-zero. Define the functional
\begin{eqnarray}\label{Phi}
\Phi_\eps(\tau,h)&=&\sum_{k=1}^\infty
h_k\bigg(\int_{-1/2}^{1/2}\cos[2\pi k(t-\tau)]\,dx^\eps(t)\bigg)^2.
\end{eqnarray}
The PMLE  of $\theta$ is then
$\hat\theta_{\text{\tiny PMLE}}=\arg\max_{\tau}\Phi_\eps(\tau,h)$.
The role of the sequence $h$ is thus to filter out the irrelevant
terms in the right side of (\ref{ML}), that is to assign a value
$h_k$ close to zero to the terms corresponding to a small
signal-to-noise ratio $|f_k|/\eps$.

For deterministic filters $h$, the asymptotic behavior of the
estimator $\hat\theta_{\text{\tiny PMLE}}$ is studied in \cite{DGT}. Under some
smoothness assumptions on $f$, for a broad choice of filters $h$,
$\hat\theta_{\text{\tiny PMLE}}$ is proved to be first-order
asymptotic efficient. Moreover, it is shown that the second-order
term of its risk expansion is $\eps^2R^\eps[f,h]/\|f'\|^4$, where
\begin{eqnarray*}
R^\eps[f,h]&=&\sum_{k=1}^\infty (2\pi
k)^2\big[(1-h_k)^2f_k^2+\eps^2h_k^2\big].
\end{eqnarray*}
This result suggests to use the filter $h_{opt}=\arg\min_h
R^\eps[f,h]$ for defining the PMLE of $\theta$. However, this
minimizer is inapplicable since it depends on $f$. To get rid of
this dependence, the minimax approach recommends the utilization of
the filter $h_{\mathcal F}=\arg\min_h \sup_{f\in \mathcal F}
R^\eps[f,h]$. If $\mathcal F$ is a ball in a Sobolev space, a
solution of this minimization problem is given by the Pinsker filter
\cite{PIN}. Although the latter leads to a second-order
minimax estimator of $\theta$ \cite[Thm.\ 2 and 3]{DGT}, it suffers
from the well known drawbacks of the minimax theory: the obtained
estimator is pessimistic and requires the precise
knowledge of the smoothness of the unknown function.

The aim of the present paper is to propose a data-driven filter
$\hat h$ so that the resulting PMLE of $\theta$ circumvents these
drawbacks. To be more precise, let us call \textit{oracle} related to
the class of filters $\widehat{\mathcal H}$ the value
$\arg\min_{h\in \widehat{\mathcal H}} R^\eps[f,h]$. We cannot use this oracle in our estimating procedure
because it depends on the unknown $f$. Nevertheless, the estimator we propose mimics well the behavior of the
oracle related to the classe of monotone filters (cf.\ Section \ref{ss4.2}) for a fixed function $f$ and
is second-order sharp adaptive over a broad scale of Sobolev balls.

\section{PMLE based on a data-driven filter}\label{S3}

\subsection{Local properties of\/ $\Phi_\eps(\tau,h)$} Let us
introduce some auxiliary notation:
\begin{align*}
y_k(\tau)&=x_k(\tau)+ix_k^*(\tau)=\sqrt{2}\int_{-1/2}^{1/2} e^{2i\pi k(t-\tau)}dx^\eps(t),\\
z_k(\tau)&=\xi_k(\tau)+i\xi_k^*(\tau)=\sqrt{2}\int_{-1/2}^{1/2} e^{2i\pi k(t-\tau)}dW(t).
\end{align*}
We will write $x_k=x_k(0)$, $x_k^*=x_k^*(0)$ and similarly for
$\xi_k$ and $\xi_k^*$. The symmetry of $f$ implies that
$x_k(\theta)=f_k+\eps \xi_k(\theta)$ and
$x_k^*(\theta)=\eps\xi_k^*(\theta)$. Moreover, for all $\tau\in\RR$,
the random variables $\{\xi_i(\tau),\xi_i^*(\tau);\,i\in\NN\}$ are
i.i.d.\ standard Gaussian. Using this notation, the functional
$\Phi_\eps$ can be rewritten as follows
$$
\Phi_\eps(\tau,h)=\frac12\sum_{k=1}^\infty
h_k\bigg(x_k(\theta)\cos[2\pi k(\tau-\theta)]+x_k^*(\theta)\sin[2\pi
k(\tau-\theta)]\bigg)^2.
$$
Our aim is to show that under some assumptions on $h$, the function
$\Phi_\eps(\cdot,h)$ has, with a probability close to one, a local
maximum in a neighborhood of $\theta$. Note that the derivative of
the function $\tau\mapsto\Phi_\eps(\tau,h)$ is given by
\begin{align}
\Phi_\eps'(\tau,h)=-&\sum_{k=1}^\infty h_k\pi k(x_k(\theta)^2-x_k^*(\theta)^2)\sin[4\pi k(\tau-\theta)]\nonumber\\
+&\sum_{k=1}^\infty 2h_k\pi k x_k(\theta)x_k^*(\theta)\cos[4\pi k(\tau-\theta)].\label{phi'1}
\end{align}

\begin{prop}\label{prop1}
Let $\|f^{(\beta_*)}\|^2=\sum_k (2\pi k)^{2\beta_*}f_k^2<\infty$ for
some $\beta_*>1$ and set $\hat\beta_*=\beta_*\wedge 1.5$,
$N_\eps=[(\eps^2\log\eps^{-5})^{-\frac{1}{2\beta_*+1}}]$. Let
$h\in[0,1]^{N_\eps}$ be a random vector depending on $(x_k,x_k^*)$ only via
$x_k^2+{x_k^*}^2$. For any $\eps<1$, there exists an event $\mathcal
A_\eps$ such that $\Pb_{\theta,f}(\mathcal A_\eps^c)\le 2\eps^4$ and
on $\mathcal A_\eps\cap\big\{\sum_{k=1}^\infty (1-h_k)(2\pi
k)^2f_k^2\leq \|f'\|^2/4\big\}$, for all $\eps$ verifying
\begin{equation}\label{eps0}
N_\eps^{1-\hat\beta_*}\le \frac{\|f'\|}{2\|f^{(\beta_*)}\|+6\pi},
\end{equation}
the function $\tau\mapsto\Phi_\eps(\tau,h)$ is strictly concave and
admits a unique maximum $\hat\theta_\eps$ in the interval
$[\theta-(4\pi N_\eps)^{-1},\theta+(4\pi N_\eps)^{-1}]$, satisfying
\begin{equation}\label{cons}
|\hat\theta_\eps-\theta|\leq 8\eps\sqrt{\log\eps^{-5}}\;\|f'\|^{-1}.
\end{equation}
\end{prop}
\begin{proof}
Set $\Theta_\eps=[\theta-(4\pi N_\eps)^{-1},\theta+(4\pi
N_\eps)^{-1}]$. Assume that (\ref{eps0}) is fulfilled and
$\sum_{k=1}^\infty (1-h_k)(2\pi k)^2f_k^2\leq \|f'\|^2/4$. On the
one hand, the first inequality of Proposition~\ref{prop2} (see Section~\ref{s6.2} below) implies
that there exists an event $\mathcal A_\eps$, such that on this event
$$
\Phi_\eps''(\tau,h)\le -\|f'\|^2+\sum_{k=1}^\infty (1-h_k)(2\pi
k)^2f_k^2+\|f'\|^2/4\le -\|f'\|^2/2,
$$
for all $\tau\in\Theta_\eps$ and for $\eps$ small enough. Therefore,
$\Phi_\eps(\cdot,h)$ is strictly concave. On the other hand, the
second inequality of Proposition~\ref{prop2} implies that
\begin{align}
\frac{\Phi_\eps'(\tau,h)}{\tau-\theta}&\le
-\frac{\|f'\|^2}{2}+\frac{2\eps(\|f'\|\sqrt{\log\eps^{-5}}+2\pi
N_\eps^{1-\hat\beta_*})}{|\tau-\theta|}\nonumber\\
&\le
-\frac{\|f'\|^2}{2}+\frac{4\eps\sqrt{\log\eps^{-5}}\;\|f'\|}{|\tau-\theta|}\label{xxx}
\end{align}
for sufficiently small values of $\eps$. Therefore,
$\pm\Phi_\eps'(\theta\pm (4\pi N_\eps)^{-1},h)<0$, which guarantees
that the maximum $\hat\theta_\eps$ of $\Phi_\eps(\cdot,h)$ is
attained in the interior of $\Theta_\eps$ and
$\Phi_\eps'(\hat\theta_\eps,h)=0$. Applying (\ref{xxx}) to
$\tau=\hat\theta_\eps$ we get (\ref{cons}).
\end{proof}
\begin{remark}\label{remN}
The choice $N_\eps=[(\eps^2\log\eps^{-5})^{-\frac{1}{2\beta_*+1}}]$
has a simple interpretation. If $f_k^2\le \eps^2$ for some
$k\in\NN$, then the $k$th observation in (\ref{x_k}) is not relevant
for estimating the parameter $\theta$. Let $\mathcal K=
\{k\in\NN:f_k^2>\eps^2\}$ and $K=\#\mathcal K$. Then
$$
\sum_{k\in\mathcal K}(2\pi k)^{2\beta_*}f_k^2\geq \eps^2\sum_{k\in
\mathcal K}(2\pi k)^{2\beta_*}\geq \eps^2\sum_{k=1}^K(2\pi
k)^{2\beta_*}\geq\frac{\eps^2(2\pi)^{2\beta_*}K^{2\beta_*+1}}{2\beta_*+1}.
$$
Thus, the number $K$ of Fourier coefficients $f_k$ larger than
$\eps$ is at most $O(\eps^{-\frac{2}{2\beta_*+1}})$. Thus, for
$\eps$ small enough, all observations relevant for estimating
the parameter $\theta$ lie in $\{y_1,\ldots,y_N\}$.
\end{remark}

\begin{remark}
In \cite{DGT}, the estimator $\hat\theta_{\text{\tiny PMLE}}$ is defined as the
maximizer of $\Phi_\eps(\cdot,h)$ over the whole interval $\Theta$.
Instead, we define it as the local maximizer in the neighborhood of
a preliminary estimator. This modification is explained by the fact
that the function $\tau\mapsto \Phi_\eps(\tau,h)$ is only locally
concave when $\beta_*<2$. Furthermore, the computation of the local
minimum is faster than the computation of the global minimum.
\end{remark}

\begin{remark}
Elaborating on the arguments of Proposition~\ref{prop1}, it can be
shown that Thm.\ 1 from \cite{DGT} remains true for $\beta>1$, provided
that $\theta_{\text{\tiny PMLE}}$ is defined as the local maximizer of
$\Phi_\eps(\cdot,h)$ and condition B2 is replaced by $h_k=0$ for any
$k>N_\eps$.
\end{remark}
\subsection{Blockwise constant Stein filter with penalization}

Let $J$ be a positive integer and
$\kappa_1,\ldots,\kappa_{J+1}\in\{1,\ldots,N_\eps\}$ be a strictly
increasing sequence such that $\kappa_1=1$. Set $B_j=\{k\in\NN:\,
\kappa_j\le k<\kappa_{j+1}\}$. Let $\mathcal H^*(B)$ be the set of
all filters $h\in[0,1]^{N_\eps}$ that are constant on the blocks
$B=\{B_j\}_{j=1}^J$:
$$
h\in \mathcal H^*(B)\qquad\Longleftrightarrow\qquad h_k=h_{k'},\quad\forall\ k,k'\in B_j.
$$
An oracle knowing the function $f$, would choose the best possible filter
$h^*\in\mathcal H^*(B)$ by minimizing $R^\eps[f,h]$ over $\mathcal
H^*(B)$. We call $h^*$ oracle choice of filter or simply oracle. Simple computations show that
$$
h_k^*=\frac{\|f'\|_{(j)}^2}{\|f'\|_{(j)}^2+\eps^2\sigma_j^2},\qquad k\in B_j,\quad j=1,\ldots,J,
$$
where $\|f'\|_{(j)}^2=\sum_{k\in B_j} (2\pi k)^2f_k^2$ and $\sigma_j^2=\sum_{k\in B_j}(2\pi k)^2$.

Since the oracle $h^*$ depends on $f$, we replace it by a suitable
estimator. Let us define
$R^\eps_j[f,a]=(1-a)^2\|f'\|^2_{(j)}+\eps^2a^2\sigma_j^2$ so that
$R^\eps[f,h]=\sum_jR^\eps_j[f,h_{\kappa_j}]+\sum_{k>N_\eps}(2\pi
k)^2f_k^2$. Then $h_{\kappa_j}^*$ is the minimizer of
$R^\eps_j[f,a]$ over $a\in\RR$. Define $\|y'\|_{(j)}^2=\sum_{k\in B_j}(2\pi k)^2|y_k|^2$ and
$y_k=x_k+ix_k^*$. For small values of
$h^*_{\kappa_j}$, the minimizer of
$$
\hat R^\eps_j[a]:=(1-a)^2(\|y'\|^2_{(j)}-2\eps^2\sigma_j^2)_++\eps^2a^2\sigma_j^2,
$$
which is an estimator of $R^\eps_j[f,a]$, can be large with respect
to $h_{\kappa_j}^*$. To avoid such a configuration, we penalize large
values of $a$ and define the estimator $\hat h^{JS}$ of $h^*$ as the
minimum over $[0,1]$ of the function:
$$
a\mapsto \hat R^\eps_j[a]+2\varphi_j\eps^2\sigma_j^2 a,
$$
where $\varphi_j>0$ is a factor of penalization tending to zero as $\eps\to 0$.
This leads us to the penalized Stein filter
\begin{align}\label{hJS}
\hat h_k^{JS}&=\bigg(1-\frac{\eps^2\sigma_j^2(1+\varphi_j)}{(\|y'\|_{(j)}^2-2\eps^2\sigma_j^2)_++\eps^2\sigma_j^2}\bigg)_+
,\quad \forall k\in B_j.
\end{align}

\subsection{Weakly geometrically increasing blocks}
The aim of this section is to propose a concrete scheme for defining
the blockwise constant data-driven filter. We use the weakly
geometrically increasing blocks introduced by Cavalier and Tsybakov
\cite{CavTsy2,CavTsy1}. These blocks have the advantage of being simple (the
construction is driven by only one parameter) and of having good
approximation properties with respect to the class of monotone
filters \cite[Lemma 1]{CavTsy1}.

Let $\nu=\nu_\eps$ be a positive integer that increases as $\eps$
decreases. Set $\rho_\eps=\nu_\eps^{-1/3}$ and define
\begin{equation}\label{kappaj}
\kappa_j=
\begin{cases}
(1+\nu_\eps)^{j-1},\ &j=1,2,\\
\kappa_{j-1}+\lfloor\nu_\eps\rho_\eps(1+\rho_\eps)^{j-2}\rfloor,\
&j=3,4,\ldots,
\end{cases}
\end{equation}
where $\lfloor x\rfloor$ stands for the largest integer strictly
smaller than $x$. Let $J$ be the smallest integer $j$ such that
$\kappa_j\ge N_\eps+1$. We redefine $\kappa_{J+1}=N_\eps+1$ and set
$B_j=\{\kappa_j,\ldots,\kappa_{j+1}-1\}$ for all $j=1,\ldots,J$.

\subsection{Brief description of the procedure}\label{Brief}
The outlined scheme can be implemented as follows.
\begin{itemize}
\item[1.] Choose a real number $\beta_*>1$ and set
$N_\eps=5\vee [(\eps^2\log\eps^{-5})^{-\frac{1}{2\beta_*+1}}]$,
$\nu_\eps=[e^{\sqrt{\log N_\eps}}]$ and $\rho_\eps=\nu_\eps^{-1/3}$.
\item[2.] Define the sequence $(\kappa_j)_j$ by (\ref{kappaj}).
\item[3.] Set $\varphi_j=\sqrt{24\log\eps^{-5}/(\kappa_{j+1}-\kappa_j)}$,
$\sigma_j^2=\sum_{\kappa_j\le k<\kappa_{j+1}} (2\pi k)^2$ and define
the data-dependent filter $\hat h^{JS}$ by (\ref{hJS}).
\item[4.] Compute the preliminary estimator $\bar\theta_\eps$ (cf.\ Section~\ref{S6}) and set
$\bar\Theta_\eps=[\bar\theta_\eps-\delta_\eps,\bar\theta_\eps+\delta_\eps]$
with $\delta_\eps=\eps\log(\eps^{-2})$.

\item[5.] Define $\hat\theta_\eps$ as the minimum in $\bar\Theta_\eps$ of $\Phi_\eps(\cdot,\hat
h^{JS})$ (see (\ref{Phi})).
\end{itemize}
Note that the only ``free'' parameter in this procedure is
$\beta_*$. In practice, if no information on the regularity of $f$
is available, it appears plausible to assume that $f$ has Sobolev
smoothness $\beta_*=2$.

\section{Main results}\label{S4}

\subsection{Comparison with the blockwise constant oracle}

In this section, $\hat h^{JS}$ denotes the blockwise constant filter
defined by (\ref{hJS}), $\varphi_j$ is the penalization we use on
the block $B_j$ and $\varphi_\eps=\max_j \varphi_j$. We emphasize
that in this section no condition on the blocks $B_j$ is required.
Let $T_j$ be the length of the block $B_j$ and $T_\eps=\inf_j T_j$.
The oracle choice of $h$ in the class $\mathcal H^*(B)$ of all
filters constant on the blocks $B=\{B_j\}_j$ is denoted by $h^*$.
Define
\begin{equation}\label{thetaJS}
\hat\theta_\eps^{JS}=\argmax_{\tau\in\bar\Theta_\eps}
\Phi_\eps(\tau,\hat h^{JS}),
\end{equation}
where $\bar\Theta_\eps=[\bar\theta_\eps-\delta_\eps,\bar\theta_\eps+\delta_\eps]$
and $\bar\theta_\eps$ is a rate optimal initial estimator of
$\theta$ (cf.\ Section~\ref{InEst}).
Introduce the functional class
$$
\mathcal F(\beta_*,L_*,\rho)=\bigg\{f\in\mathcal F_0(\rho):\
\|f^{(\beta_*)}\|\leq L_*  \bigg\},
$$
where $\beta_*>1$, $\rho>0$, $L_*>0$ are some constants.
\begin{theorem}\label{Th1}
Let  $\hat\theta_\eps^{JS}$ be defined by (\ref{thetaJS}) with
blocks $B_j$ verifying $\log\eps^{-1}=o(T_\eps)$ as
$\eps\to 0$. If the penalty $\varphi_j$ is equal to
$\sqrt{24T_j^{-1}\log\eps^{-5}}$, then
$$
\eps^{-2}\|f'\|^2\esp_{\theta,f}[(\hat\theta_\eps^{JS}-\theta)^2]\leq
1+(1+\alpha_\eps)\frac{R^\eps[f,h^*]}{\|f'\|^2},
$$
where $\alpha_\eps\to 0$ as $\eps\to 0$ uniformly in $f\in \mathcal
F(\beta_*,L_*,\rho)$.
\end{theorem}

\begin{remark}
If the block $B_j$ is large, then more observations $(x_k,x_k^*)$
are used for estimating the value of the oracle $h^*_{\kappa_j}$.
Hence, it is natural to expect that $\alpha_\eps$ decreases as
$T_\eps$ increases. A thorough inspection of the proof allows to
describe this feature with the help of the order
relation $\alpha_\eps^2\asymp{T_\eps^{-1}\log\eps^{-1}}$.
\end{remark}

\begin{proof} Let us denote
$\mathcal E=\{\hat
h_k^{JS}\in[h_k^*-\sqrt8\varphi_{j(k)}(1-h_k^*),h_k^*],\ \forall
k=1,\ldots, N_\eps\}$ and $\tilde h=\hat h^{JS}\1_{\mathcal E}$,
where $j(k)$ is the number of the block containing $k$.
Lemma~\ref{lemh_k} implies that $\Pb_{\theta,f}(\mathcal E^c)\le
2\eps^4$. For $\eps$ small enough, we have $(1-\tilde h_k)\1_{\mathcal E}\leq
2(1-h_k^*)$ and the inequality $\sum_k(1-\hat h_k)(2\pi k)^2f_k^2\le
\|f'\|^2/4$ is fulfilled on $\mathcal E$.

By virtue of Proposition~\ref{prop3}, for sufficiently small values
of $\eps$, the event $\mathcal A_0=\Big\{|\bar\theta_\eps-\theta|\le
\big(\frac1{4\pi
N_\eps}-\delta_\eps\big)\wedge\big(\delta_\eps-\frac{8\eps\sqrt{\log\eps^{-5}}}{\|f'\|}\big)\Big\}$
verifies $\Pb_{\theta,f}(\mathcal A_0^c)=O(\eps^4)$. It can be
checked that on $\mathcal A_0$,
$
\Big[\theta\pm \frac{8\eps\sqrt{\log\eps^{-5}}}{\|f'\|}\Big]
\subset [\bar\theta_\eps\pm\delta_\eps]\subset
\big[\theta\pm\frac1{4\pi N_\eps}\big],
$
where we have used the notation $[\theta\pm
\delta]:=[\theta-\delta,\theta+\delta]$.

According to Proposition~\ref{prop1}, on the event $\mathcal
A_1=\mathcal A_0\cap \mathcal A_\eps\cap \mathcal E$, the function
$\Phi_\eps(\cdot,\hat h)$ is strictly concave and has a unique
maximum in $[\bar\theta_\eps\pm\delta_\eps]$. Hence the estimator
$\hat\theta_\eps^{JS}$ verifies
$\Phi_\eps'(\hat\theta_\eps^{JS},\tilde h)=0$ on the event $\mathcal
A_1$. By Taylor's formula, there exists
a point $\tilde\theta\in[\theta,\hat\theta_\eps^{JS}]$ such that
$$
0=\Phi_\eps'(\hat\theta_\eps^{JS},\tilde h)=\Phi_\eps'(\theta,\tilde
h)+(\hat\theta_\eps^{JS}-\theta)\Phi_\eps''(\theta,\tilde h)+
\frac{(\hat\theta_\eps^{JS}-\theta)^2}{2}\Phi_\eps'''(\tilde\theta,\tilde
h).
$$
Therefore, on $\mathcal A_1$,
\begin{align*}
\hat\theta_\eps^{JS}-\theta&=-\frac{\Phi_\eps'(\theta,\tilde h)}
{\Phi_\eps''(\theta,\tilde
h)+\frac12(\hat\theta_\eps^{JS}-\theta)\Phi_\eps'''(\tilde\theta,\tilde
h)}\ .
\end{align*}
Using (\ref{phi'1}), one checks that
\begin{align*}
\Phi''_\eps(\theta,\tilde h)&=-\sum_{k=1}^\infty \tilde h_k(2\pi
k)^2[(f_k+\eps\xi_k(\theta))^2-\eps^2\xi_k^*(\theta)^2].
\end{align*}
In the sequel, we write $\xi_k,\xi_k^*$ instead of $\xi_k(\theta),\xi_k^*(\theta)$.
On the one hand, Lemmas \ref{lem10}, \ref{lem11} (with $x^2=\log\eps^{-5}$) and \ref{lem12}
combined with (\ref{1-h^*}) imply that, on an event $\mathcal A_2$
of probability higher than $1-4\eps^4$, we have
\begin{align*}
\Phi_\eps''(\theta,\tilde h)&=
-\|f'\|^2+\sum_{k=1}^\infty(1-h_k^*)(2\pi k)^2f_k^2-2\eps\sum_{k=1}^{N_\eps} h_k^*(2\pi k)^2f_k\xi_k+o(R^\eps)\\
&=-\|f'\|^2(1-\|f'\|^{-2}R^\eps[f,h^*]-\zeta-o(R^\eps[f,h^*])),
\end{align*}
where $\zeta=2\eps\|f'\|^{-2}\sum_k h_k^*(2\pi k)^2 f_k\xi_k$ is a zero mean Gaussian random variable.
By virtue of (\ref{3Tj}), its variance verifies
\begin{equation}\label{zeta}
4\eps^2\sum_{k=1}^{N_\eps} {h_k^*}^2(2\pi k)^4f_k^2\leq
\frac{12\eps^2}{\min_jT_j}\sum_{j=1}^J\hkapjsq\sigma_j^2\|f'\|_{(j)}^2\leq
\frac{12\|f'\|^2R^\eps[f,h^*]}{T_\eps}.
\end{equation}
Therefore, by Rosenthal's inequality, $\esp_{\theta,f}[\zeta^{2p}]^{1/p}=o(R^\eps[f,h^*])$ for
any $p>0$. On the other hand, in view of Lemma~\ref{lemPhi'''},
there is an event $\mathcal A_3$ such that $\Pb_{\theta,f}(\mathcal
A_3^c)=O(\eps^4)$ and
$$
(\hat\theta_\eps^{JS}-\theta)\Phi_\eps'''(\tilde\theta,\tilde
h)=o(R^\eps)
$$
on this event. Using the inequality $(1-x)^{-2}\leq 1+2x+16x^2$ for
all $x\in [-1/2,1/2]$, we get
\begin{align*}
\|f'\|^4(\hat\theta_\eps^{JS}-\theta)^2&=\frac{\Phi_\eps'(\theta,\tilde
h)^2}
{\big(1-\|f'\|^{-2}R^\eps[f,h^*]-\zeta-o(R^\eps)\big)^2}\\
&\leq \Phi_\eps'(\theta,\tilde h)^2\bigg(1 +\frac{2R^\eps[f,h^*]}{\|f'\|^{2}}+2\zeta+o(R^\eps)\bigg)+o(\eps^2R^\eps),
\end{align*}
on the event $\mathcal A_4=\mathcal A_1\cap\mathcal A_2\cap\mathcal A_3$.
Using Lemma~\ref{lemPhi'}, we infer that
\begin{align*}
\esp_{\theta,f}[\Phi'_\eps(\theta,\tilde h)^2]&\leq \eps^2\|f'\|^2+\eps^2\sum_{k=1}^\infty(2\pi k)^2[({h_k^*}^2-1)f_k^2+{h_k^*}^2\eps^2]\\
&=\eps^2(\|f'\|^2-R^\eps[f,h^*]).
\end{align*}
Combining these relations with Lemmas~\ref{lem13} and
\ref{lemPhi'z}, we get an event $\mathcal A$ such that
$\Pb_{\theta,f}(\mathcal A^c)=O(\eps^4)$ and
\begin{align*}
\|f'\|^4\esp_{\theta,f}[(\hat\theta_\eps^{JS}-\theta)^2\1_{\mathcal
A}]&\leq \esp_{\theta,f}[\Phi'_\eps(\theta,\tilde h)^2]
\bigg(1+\frac{2R^\eps[f,h^*]}{\|f'\|^2}\bigg)\!\!+o(\eps^2R^\eps[f,h^*])\\
&\leq \eps^2 \|f'\|^2+\eps^2R^\eps[f,h^*]+o(\eps^2R^\eps[f,h^*]).
\end{align*}
Since $|\hat\theta^{JS}_\eps-\theta|\le
1$, we have $\esp_{\theta,f}[(\hat\theta_\eps^{JS}-\theta)^2\1_{\mathcal A^c}]\le
\Pb_{\theta,f}(\mathcal A^c)= O(\eps^{4})$.
In view of (\ref{R^e}), $\eps^2=o(R^\eps[f,h^*])$. Therefore
$\esp_{\theta,f}[(\hat\theta_\eps^{JS}-\theta)^2\1_{\mathcal
A^c}]=O(\eps^4)=o(\eps^2R^\eps[f,h^*])$ and the assertion of the
theorem follows.
\end{proof}

\subsection{Comparison with the monotone oracle}\label{ss4.2}

Now we consider the class $\mathcal H_{\rm mon}$ of filters having
decreasing components, that is
$$
\mathcal H_{\rm mon}=\Big\{h\in[0,1]^{N_\eps}\;:\; h_k\ge h_{k+1},\ 1\le
k\le N_\eps-1 \Big\}.
$$
The class
$\mathcal H_{\rm mon}$ is of high interest in statistics because it
contains the most common filters such as the projection filter, the
Pinsker filter, the Tikhonov or smoothing spline filter and so
forth.

\begin{prop}\label{prop4}
Set $\gamma_\eps=\max_{1\le j\le J-1}
({\sigma_{j+1}^2}/{\sigma_j^2})$. Then
$$
R^\eps[f,h^*]\leq \gamma_\eps\inf_{h\in\mathcal H_{\rm mon}}
R^\eps[f,h]+\eps^2\sigma_1^2.
$$
\end{prop}

A more general version of this result is Lemma~1  in \cite{CavTsy1}.
Since the proof in our setting is simple, we give
it below.

\begin{proof}
Let $h$ be a filter from $\mathcal H_{\rm mon}$. Define $\bar h$ by
$\bar h_k=h_{\kappa_j}$ if $k\in B_j$ for some $j$ and $\bar h_k=0$
if $k>N_\eps$. Since the components of $h$ are decreasing, we have
$1-\bar h_k\leq 1-h_k$ and therefore
\begin{align}\label{barh}
R^\eps[f,\bar h]\leq \sum_{k=1}^{N_\eps}(1- h_k)^2(2\pi
k)^2f_k^2+\eps^2\sum_{j=1}^J h_{\kappa_j}^2\sigma_j^2.
\end{align}
Again by monotonicity of $h$, we have $h_{\kappa_j}\leq h_k$ for all
$k\in B_{j-1}$. Hence,
$$
\sum_{j=2}^J h_{\kappa_j}^2\sigma_j^2\leq \gamma_\eps\sum_{j=2}^J
h_{\kappa_j}^2\sigma_{j-1}^2\leq \gamma_\eps\sum_{k=1}^{N_\eps} h_k^2(2\pi
k)^2.
$$
Combining this inequality with (\ref{barh}) and bounding $h_1$ by
$1$, we get $R^\eps[f,\bar h]\leq \gamma_\eps
R^\eps[f,h]+\sigma_1^2$. Since $\bar h\in \mathcal H^*$ and $h^*$
minimizes $R^\eps[f,h]$ over all $h\in \mathcal H^*$, we have
$R^\eps[f,h^*]\leq \gamma_\eps R^\eps[f,h]+\eps^2\sigma_1^2$. This
inequality holds for every $h\in\mathcal H_{\rm mon}$, therefore
the assertion of the proposition follows.
\end{proof}

Combining this proposition with Theorem~\ref{Th1} we get the
following result.

\begin{corollary}\label{cor2}
Assume that the conditions of Theorem~\ref{Th1} are fulfilled, then
$$
\eps^{-2}\|f'\|^2\esp_{\theta,f}[(\hat\theta_\eps^{JS}-\theta)^2]\leq
1+\gamma_\eps(1+\alpha_\eps)\frac{\min_{h\in \mathcal H_{\rm mon}}
R^\eps[f,h]}{\|f'\|^2},
$$
where $\alpha_\eps\to 0$ as $\eps\to 0$ uniformly in $f\in \mathcal
F(\beta_*,L_*,\rho)$.
\end{corollary}

\begin{remark}
For the blocks defined by (\ref{kappaj}), we have
$T_\eps=\nu_\eps\rho_\eps(1+\rho_\eps)$, $\sigma_1^2\leq
4\pi^2\nu_\eps^3$ and
$-\nu_\eps\rho_\eps+\nu_\eps(1+\rho_\eps)^j\le\kappa_{j+1}\le
1+\nu_\eps(1+\rho_\eps)^j$. One also checks that $\gamma_\eps=\max_j
\sigma_{j+1}^2/\sigma_j^2$ is asymptotically equivalent to
$(1+\rho_\eps)^3\sim 1+3\rho_\eps$ as $\eps\to 0$. Therefore the
factor in the oracle inequality of Corollary~\ref{cor2} is of order
$(1+3\rho_\eps+\alpha_\eps)$. We have already mentioned that
$\alpha_\eps^2=O(T_\eps^{-1}\log\eps^{-1})$. The trade-off between $\alpha_\eps$ and
$\rho_\eps$ leads us to $\rho_\eps\asymp \nu_\eps^{-1/3}$.
This clarifies our choice of $\rho_\eps$ slightly differing from the one of
\cite{CavTsy1}.
\end{remark}
\begin{remark}
In \cite{CavTsy2,RIG, TSY} the weakly geometrically
increasing blocks are defined by $T_j=\lfloor\nu(1+\rho)^{j-1}\rfloor$. This
type of blocks does not lead to a sharp oracle inequality in our case, since
we need not only $\max(T_{j+1}/T_j)\to 1$, but also $\max(\kappa_{j+1}/\kappa_j)\to 1$
as $\eps\to 0$.
\end{remark}

\subsection{Second-order minimax sharp adaptation}

To complete the theoretical analysis, we show below that the
estimator $\hat\theta_\eps^{JS}$ corresponding to the blocks
(\ref{kappaj}) enjoys minimax properties over a large scale of
Sobolev balls. Assume that $\bar f\in \mathcal F(\beta^*,L^*,\rho)$
and define
$$
\mathcal F_{\delta,\beta,L}(\bar f)=\Big\{f=\bar f+v\,:\;\|v\|\leq \delta,\
\|v^{(\beta)}\|\leq L\Big\}.
$$
\begin{theorem}\label{Th2}
Assume that $\nu_\eps$ verifies $\eps^\frac2{2\beta+1}\nu_\eps\to 0$
as $\eps\to 0$ and the conditions of Theorem~\ref{Th1} are
fulfilled. If $\delta=\delta_\eps$ tends to zero as $\eps\to 0$ and $\bar f\in\mathcal
F(\beta^*,L^*,\rho)$ with $\beta^*>\beta\geq \beta_*$, then the
estimator $\hat\theta_\eps^{JS}$ defined in Section~\ref{Brief} satisfies
$$
\sup_{\theta\in\Theta,f\in \mathcal F_{\delta,\beta,L}(\bar f)}
\eps^{-2}\|f'\|^2\esp_{\theta,f}[(\hat\theta_\eps^{JS}-\theta)^2]\leq
1+(1+o(1))\frac{C(\beta,L)\eps^{\frac{4\beta-4}{2\beta+1}}}{\|\bar
f'\|^2},
$$
when $\eps\to 0$, with $ C(\beta,L)=\frac13\big(\frac{\beta-1}{2\pi
(\beta+2)}\big)^\frac{2\beta-2}{2\beta+1}(L(2\beta+1))^\frac3{2\beta+1}$.
Moreover, the following lower bound holds:
$$
\inf_{\tilde\theta_\eps}\sup_{\theta\in\Theta,f\in \mathcal F_{\delta,\beta,L}
(\bar f)}
\eps^{-2}\|f'\|^2\esp_{\theta,f}[(\tilde\theta_\eps-\theta)^2]\geq
1+(1+o(1))\frac{C(\beta,L)\eps^{\frac{4\beta-4}{2\beta+1}}}{\|{\bar
f}'\|^2},
$$
where the $\inf$ is taken over all possible estimators
$\tilde\theta_\eps$.
\end{theorem}

\begin{proof}
According to Lemma~\ref{lemPinsk}, there exists a
filter $\lambda^*\in\mathcal H_{\rm mon}$ such that
$$\sup_{f\in
\mathcal F_{\delta,\beta,L}(\bar f)}
R^\eps[f,\lambda^*]/\|f'\|^2\leq (1+o(1))C(\beta,L)\eps^{\frac{4\beta-4}{2\beta+1}}\sup_f\|
f'\|^{-2}.
$$
Since $\delta_\eps\to 0$ as $\eps\to 0$, we have $\sup_{f\in\mathcal F_{\delta,\beta,L}(\bar f)}\|f'\|^{-2}=(1+o(1))\|\bar f'\|^{-2}$.
Combining this result with Corollary~\ref{cor2}, we get
the first inequality.
The second inequality is Theorem 2 of \cite{DGT}. Although the latter is stated for $\beta\geq 2$, the
inspection of its proof shows that the same claim is true for any $\beta>1$.
\end{proof}

\section{Preliminary estimator and technical lemmas}\label{S6}

\subsection{Preliminary estimator}\label{InEst} Having the observation
$(x^\eps(t),\,|t|\le1/2)$, we can compute $ x_1, x_1^*$ by
(\ref{xx*}). Then we have $x_1=f_1\cos(2\pi\theta)+\eps\xi_1$ and
$x_1^*=f_1\sin(2\pi\theta)+\eps\xi_1^*$, where $\xi_1$, $\xi_1^*$ are
independent standard Gaussian random variables. We define
$$
\bar\theta_\eps=\frac1{2\pi}\;\arctan\bigg(\frac{x_1^*}{x_1}\bigg),
$$
if $x_1\not=0$ and $\bar\theta_\eps=1/4$ if $x_1=0$. One easily
checks that $\bar\theta_\eps$ is the maximum likelihood estimator in
the model induced by observations $(x_1,x_1^*)$. The following
result describes its asymptotic behavior.

\begin{prop}\label{prop3}
If $\eps$ is sufficiently small, then
\begin{align*}
\sup_{|\theta|\le T}\Pb_{\theta,f}\big(|\bar\theta_\eps-\theta|\ge x
\big)&\le \exp\big(-2(x/\eps)^2f_1^2 \cos^2(2\pi T)\big),
\end{align*}
for all $x\in[0,1/2]$ and for all $T<1/4$.
\end{prop}
\begin{proof}
Let us introduce $X=\sqrt{\xi_1^2+{\xi_1^*}^2}$. One checks that
$\sin[2\pi(\bar\theta_\eps-\theta)]=\eps(\xi_1\sin(2\pi\bar\theta_\eps)-\xi_1^*\cos(2\pi\bar\theta_\eps))f_1^{-1}$.
The Cauchy-Schwarz inequality implies that
$|\sin[2\pi(\bar\theta_\eps-\theta)]|\le \eps X|f_1^{-1}|$. Since
$\theta\in[-1/4,1/4]$ and $\bar\theta_\eps\in [-T,T]$, we have
$\frac{\sin[2\pi(\bar\theta_\eps-\theta)]}{2\pi(\bar\theta_\eps-\theta)}\geq
\cos(2\pi T)$. Therefore,
\begin{align*}
\Pb_{\theta,f}\big(|\bar\theta_\eps-\theta|\ge x
\big)&\le\Pb_{\theta,f}\bigg(\frac{|\sin[2\pi(\bar\theta_\eps-\theta)]|}{2\pi\cos(2\pi T)}\ge
x\bigg)\\
&\le\Pb_{\theta,f}\big(X\ge 2x\eps^{-1}\pi|f_1| \cos(2\pi T)\big),
\end{align*}
and the fact that $X^2/2$ follows the exponential law completes
the proof.
\end{proof}

\subsection{Proofs of Lemmas used in Proposition~\ref{prop1}}\label{s6.2}
Let us start with some basic facts that will be often used in the
proofs. For any $n,m,p\in\NN$, we have
\begin{align}\label{mnp}
n^p(n-m)\ge \sum_{k=m+1}^n k^p\ge \frac{n^p(n-m)}{p+1}.
\end{align}
Applying this inequality to $p=2$, we get
\begin{align}\label{3Tj}
\max_{k\in B_j}(2\pi k)^2\leq 3\sigma_j^2/T_j.
\end{align}
Assume now that $\xi$ is a random variable of law $\mathcal N(0,1)$.
For any $\sigma^2\le1/4$, we have $(1-2\sigma^2)^{-1}\le 2$ and
$(1-2\sigma^2)^{-1/2}\le e^{2\sigma^2}$, therefore
\begin{equation}\label{chi2}
\esp[e^{(\mu
+\sigma\xi)^2}]=\frac{e^{\frac{\mu^2}{(1-2\sigma^2)}}}{\sqrt{1-2\sigma^2}}\le
e^{\frac{\mu^2}{(1-2\sigma^2)}+2\sigma^2}\leq
\exp\big(2\mu^2+2\sigma^2\big).
\end{equation}
Using the more precise inequalities $(1-2\sigma^2)^{-1}\le
1+4\sigma^2$ and $\log(1-2\sigma^2)^{-1}\leq 2\sigma^2+4\sigma^4$,
we get $\esp[e^{(\mu+\sigma\xi)^2}]\le
\exp({\mu^2+\sigma^2+2\sigma^2(2\mu^2+\sigma^2)})$ or equivalently,
\begin{equation}\label{chi2'}
\esp[e^{2\mu\sigma\xi +\sigma^2(\xi^2-1)}]\leq
e^{2\sigma^2(2\mu^2+\sigma^2)}.
\end{equation}

Throughout this section, we assume that $\|f^{(\beta_*)}\|^2=\sum_k
(2\pi k)^{2\beta_*}f_k^2<\infty$ for some $\beta_*>1$,
$N_\eps=[(\eps^2\log\eps^{-5})^{-\frac{1}{2\beta_*+1}}]$ and
$h\in[0,1]^{N_\eps}$ is a random vector depending on $(x_k,x_k^*)$ only via
$x_k^2+{x_k^*}^2$. Without loss of generality, we give the proofs in
the case $\theta=0$.

\begin{lemma}\label{lem1} Set $\hat\beta_*=\beta_*\wedge 1.5$. For all $\tau$
such that $4\pi|\tau-\theta|\leq N_\eps^{-1}$,
$$
\bigg|\Phi_\eps''(\tau,h)+\sum_{k=1}^{N_\eps} (2\pi k)^2h_kf_k^2\bigg|\leq
N_\eps^{2-2\hat\beta_*}\Big(\|f^{(\beta_*)}\|+2\pi\hat
X/\sqrt{\log\eps^{-5}}\,\Big)^2,
$$
where $\hat X=\max_{1\le k\le N_\eps}\sqrt{\xi_k^2+{\xi_k^*}^2}$.
\end{lemma}

\begin{proof}
One easily checks that
\begin{align*}
\Phi_\eps''(\theta+\tau,h)=\,&-\sum_{k=1}^{N_\eps}(2\pi k)^2
h_kf_k^2\cos[4\pi k\tau]\\
&-2\eps\sum_{k=1}^{N_\eps}(2\pi
k)^2h_kf_k(\xi_k\cos[4\pi k\tau]+\xi_k^*\sin[4\pi k\tau])\\
&-\eps^2\sum_{k=1}^{N_\eps}(2\pi k)^2h_k\big[(\xi_k^2-{\xi_k^*}^2)\cos[4\pi
k\tau]+2\xi_k^*\xi_k\sin[4\pi k\tau]\big].
\end{align*}
On the one hand, thanks to inequality $|1-\cos x|\le |x|$,
\begin{align*}
\bigg|\sum_{k=1}^{N_\eps}(2\pi k)^2 h_kf_k^2(1-\cos[4\pi k\tau])\bigg|&\le
2\tau
\sum_{k=1}^{N_\eps}(2\pi k)^{3}f_k^2\\
&\le 2\tau (2\pi N_\eps)^{3-2\hat\beta_*}\|f^{(\beta_*)}\|^2\\
&\le  N_\eps^{2-2\hat\beta_*}\|f^{(\beta_*)}\|^2.
\end{align*}
On the other hand, in view of the Cauchy-Schwarz inequality, it holds
$\xi_k\cos[4\pi k\tau]+\xi_k^*\sin[4\pi k\tau]\leq \hat X$ and
$(\xi_k^2-{\xi_k^*}^2)\cos[4\pi k\tau]+2\xi_k^*\xi_k\sin[4\pi
k\tau]\leq \hat X^2$. Therefore, it holds
\begin{align*}
\bigg|\sum_{k=1}^{N_\eps}(2\pi k)^2h_kf_k&(\xi_k\cos[4\pi
k\tau]+\xi_k^*\sin[4\pi k\tau])\bigg|\le \hat X\sum_{k=1}^{N_\eps}(2\pi k)^2 |f_k| \\
&\le \hat X\|f^{(\beta_*)}\|\sqrt{\sum_{k=1}^{N_\eps}(2\pi k)^{4-2\beta_*}}
\le 2\pi\hat XN_\eps^{\frac52-\hat\beta_*} \|f^{(\beta_*)}\|,
\end{align*}
and
\begin{align*}
\bigg|\sum_{k=1}^{N_\eps}(2\pi k)^2h_k\big[(\xi_k^2-{\xi_k^*}^2)\cos[4\pi
k\tau]+2\xi_k^*\xi_k\sin[4\pi k\tau]\big]\bigg| &\le 4\pi^2
N_\eps^3\hat X^2.
\end{align*}
Taking into account the identity
$\eps^2N_\eps^3=N_\eps^{2-2\beta_*}/\log(\eps^{-5})$, for all $\tau$
verifying $|\tau|\le (4\pi N_\eps)^{-1}$, we get
\begin{align*}
\bigg|\Phi_\eps''(\theta+\tau,h)+\sum_{k=1}^{N_\eps} h_k(2\pi
k)^2f_k^2\bigg|&\leq
N_\eps^{2-2\hat\beta_*}(\|f^{(\beta_*)}\|^2+2\pi\hat X
/\sqrt{\log\eps^{-5}})^2
\end{align*}
and the assertion of the lemma follows.
\end{proof}

\begin{lemma}\label{lem2}
Let $h\in[0,1]^{N_\eps}$ be a random vector depending on $(x_k,x_k^*)$ only
via $x_k^2+{x_k^*}^2$. For any $x\in[0,\sqrt{N_\eps/6}]$, it holds
$$
\Pb_{\theta,f}\bigg(|\Phi'_\eps(\theta,h)|>2x\eps(\|f'\|+2\pi\eps
N_\eps^{3/2})\bigg)\le 2e^{-x^2}
$$
\end{lemma}

\begin{proof}
The random variables $X_k=(2\pi k)(f_k+\eps\xi_k)\xi_k^*$,
$k=1,\ldots,N_\eps$ and $h_1,\ldots,h_N$ fulfill the conditions of
Lemma~\ref{lemsym} with $\varrho_k^*=1$, $T_k=1/(2\sqrt2\pi k
\eps)$, $g_k^2=(2\pi k)^2(f_k^2+\eps^2)$, since due to (\ref{chi2}),
\begin{equation}\label{e^tX}
\esp[e^{tX_k}]=\esp[e^{(2\pi kt)^2(f_k+\eps\xi_k)^2/2}]\leq e^{(2\pi
kt)^2(f_k^2+\eps^2)}.
\end{equation}
By definition, $\Phi'_\eps(\theta,h)=\eps\sum_{k=1}^{N_\eps} h_k
X_k$, and therefore,
$$
\Pb_{\theta,f}\Big(|\Phi'_\eps(\theta,h)|\ge
2x\eps\Big(\sum_{k=1}^{N_\eps}(2\pi k)^2(f_k^2+\eps^2)\Big)^{1/2}\Big)\leq
2e^{-x^2}
$$
for all $x\in[0,(\sum_k (2\pi k)^2(f_k^2+\eps^2))^{1/2}/(2\sqrt 2
\pi N_\eps\eps)]$. To complete the proof, it suffices to remark that
$$
\frac{\sum_{k=1}^{N_\eps} (2\pi k)^2(f_k^2+\eps^2)}{8\pi^2
N_\eps^2\eps^2}\geq \frac{\sum_{k=1}^{N_\eps} k^2}{2N_\eps^2}\ge
\frac{N_\eps}6
$$
and $\sum_{k=1}^{N_\eps}(2\pi k)^2(f_k^2+\eps^2)\leq (\|f'\|+\eps(2\pi)
N_\eps^{3/2})^2$.
\end{proof}

\begin{prop}\label{prop2} Assume that $\|f^{(\beta_*}\|<\infty$ for some $\beta_*>1$
and set $\hat\beta_*=\beta_*\wedge 1.5$.
There exists an event $\mathcal A_\eps$ such that for every
$\eps<1/2$, $\Pb_{\theta,f}(\mathcal A_\eps)\ge 1-2\eps^4$ and on
$\mathcal A_\eps$ it holds:
\begin{align}
\Phi_\eps''(\tau,h)&\le -\sum_{k=1}^{N_\eps}h_k(2\pi
k)^2f_k^2+N_\eps^{2-2\hat\beta_*}\big(\|f^{(\beta_*)}\|+3\pi\big)^2\label{phi''}\\
\frac{\Phi_\eps'(\tau,h)}{\tau-\theta}&\le -\sum_{k=1}^{N_\eps}h_k(2\pi
k)^2f_k^2+N_\eps^{2-2\hat\beta_*}\big(\|f^{(\beta_*)}\|+3\pi\big)^2\nonumber\\
&\quad+\frac{2\eps(\|f'\|\sqrt{\log\eps^{-5}}+2\pi
N_\eps^{1-\beta_*})}{|\tau-\theta|}.\label{phi'}
\end{align}
for all $\tau\in[\theta-(4\pi N_\eps)^{-1},\theta+(4\pi
N_\eps)^{-1}]$.
\end{prop}
\begin{proof}
According to Lemma \ref{lem1}, we have
$$
\Phi_\eps''(\tau,h)\le -\sum_{k=1}^{N_\eps} h_k(2\pi
k)^2f_k^2+N_\eps^{2-2\hat\beta_*}\big(\|f^{(\beta_*)}\| +2\pi\hat
X/\sqrt{\log\eps^{-5}}\;\big)^2.
$$
Since for every $k$, $(\xi_k^2+{\xi_k^*}^2)/2$ follows the
exponential law with mean $1$, we have
$$
\Pb\big(4\hat X^2>9\log\eps^{-5}\big)\le N_\eps\Pb(X_1^2/2 >
\log\eps^{-5})\le N_\eps\eps^{5}\le \eps^4.
$$
This inequality completes the proof of (\ref{phi''}).

To prove (\ref{phi'}), note that for some $\tilde\tau\in[\theta,\tau]$, we have
$\Phi_\eps'(\tau,h)=\Phi_\eps'(\theta,h)+
(\tau-\theta)\Phi_\eps''(\tilde\tau,h)$.
Lemma~\ref{lem2} and (\ref{phi''}) yield (\ref{phi'}).
\end{proof}

\subsection{Lemmas used in Theorem~\ref{Th1}}

Let us start with some simple algebra allowing to obtain a rough evaluation of $R^\eps[f,h^*]$, where
$h^*$ is the ideal filter an oracle would choose in the class of blockwise constant filters.
For this filter $h^*$, it holds
\begin{align*}
R^\eps[f,h^*]=\sum_{j=1}^J\frac{\eps^2\sigma_j^2\|f'\|_{(j)}^2}{\|f'\|_{(j)}^2+\eps^2\sigma_j^2}+\sum_{k>N_\eps}
(2\pi k)^2f_k^2.
\end{align*}
Using the explicit form of $h_k^*$, we get
\begin{align}\label{1-h^*}
R^\eps[f,h^*]=\sum_{k=1}^\infty (1-h_k^*)(2\pi k)^2f_k^2\geq \eps^2\sum_{k=1}^{N_\eps} h_k^*(2\pi k)^2.
\end{align}
Since $\eps^2N_\eps^3\to 0$ as $\eps\to 0$, we have
$R^\eps[f,h^*]\leq CN_\eps^{2-2\beta_*}\to 0$ as $\eps\to 0$. In
view of $\|f'\|_{(1)}^2\ge 4\pi^2 \rho^2$ and $\eps^2\sigma_1^2\leq
4\pi^2\eps^2 N_\eps^3\to 0$ as $\eps\to 0$, for $\eps$ small enough
the inequality $\eps^2\sigma_1^2\leq \|f'\|^2_{(1)}$ holds.
Therefore,
\begin{align}\label{R^e}
R^\eps[f,h^*]\geq
\frac{\eps^2\sigma_1^2\|f'\|_{(1)}^2}{\|f'\|_{(1)}^2+\eps^2\sigma_1^2}\geq \frac{\eps^2\sigma_1^2}{2}.
\end{align}
Hence, for every function $f$, the quantity $R^\eps[f,h^*]$ tends to
zero as $\eps\to 0$ slower than $\eps^2$ and faster than
$\eps^2N_\eps^3$.

\begin{lemma}\label{lemPhi'}
It holds $\displaystyle\esp[\Phi'_\eps(\theta,\tilde h)^2]\leq
\eps^2\sum_{k=1}^{N_\eps} {h_k^*}^2(2\pi k)^2(f_k^2+\eps^2)$.
\end{lemma}
\begin{proof}
Using (\ref{phi'1}), one checks that
\begin{align}\label{hatPhi'}
\Phi'_\eps(\theta,\tilde h)&=\eps\sum_{k=1}^{N_\eps}\tilde h_k(2\pi
k)(f_k+\eps\xi_k(\theta))\xi_k^*(\theta).
\end{align}
In the sequel, we write $\xi_k,\xi_k^*$ instead of $\xi_k(\theta),\xi_k^*(\theta)$.
For any $k'\not=k$, the random variable $\tilde h_k(f_k+\eps\xi_k)\xi_k^*\tilde h_{k'}(f_{k'}+\eps\xi_{k'})\xi_{k'}^*
$ is symmetric. Therefore it has zero mean and
\begin{align*}
\esp[\Phi'_\eps(\theta,\tilde h)^2]
&=\eps^2\sum_{k=1}^{N_\eps} (2\pi k)^2\esp\big[\tilde h_k^2(f_k+\eps\xi_k)^2{\xi_k^*}^2\big]\\
&\le\eps^2\sum_{k=1}^{N_\eps} (2\pi k)^2{h_k^*}^2\esp\big[(f_k+\eps\xi_k)^2{\xi_k^*}^2\big],
\end{align*}
and the assertion of the lemma follows.
\end{proof}

\begin{lemma}\label{lemeta}
Let us denote
$$
\eta_j=\frac{2\eps\sum_{k\in B_j}(2\pi k)^2f_k\xi_k+\eps^2\sum_{k\in B_j}(2\pi k)^2(\xi_k^2+{\xi_k^*}^2-2)}
{\|f'\|_{(j)}^2+\eps^2\sigma_j^2}.
$$
For any positive $x$ such that $x^2\leq {T_j/10}$, it holds
$$
\Pb\bigg(|\eta_j|>x\sqrt{24(1-h_{\kappa_j}^*)T_j^{-1}}\bigg)\leq
2e^{-x^2}.
$$
\end{lemma}
\begin{proof}
Set $Y_k=2\eps(2\pi k)^2f_k\xi_k+\eps^2(2\pi k)^2(\xi_k^2+{\xi_k^*}^2-2)$, $\sigma=\eps(2\pi k)\sqrt{t}$ and
$\mu=(2\pi k)f_k\sqrt{t}$. Using (\ref{chi2'}) we infer that
\begin{align*}
\esp[e^{tY_k}]&=\esp[e^{2\mu\sigma\xi_k+\sigma^2(\xi_k^2-1)}]\,\esp[e^{\sigma^2({\xi_k^*}^2-1)}]\le
e^{4\sigma^2(\mu^2+\sigma^2)}\\
&=e^{4\eps^2(2\pi k)^4t^2(f_k^2+\eps^2)}
\end{align*}
as soon as $\sigma=\eps(2\pi k)\sqrt{t}\leq 1/2$, or equivalently $t\leq 1/{4\eps^2(2\pi k)^2}$. By \cite[Thm. 2.7]{PET}, for any $x>0$,
we get
$$
\Pb\bigg(\bigg|\sum_{k\in B_j}Y_j\bigg|\geq x\Big(2\sum_{k\in
B_j}4\eps^2(2\pi k)^4(f_k^2+\eps^2)\Big)^{1/2}\bigg)\leq
2e^{-x^2(1\wedge Q_\eps x^{-1})}
$$
where
$$
Q_\eps=\frac{(8\eps^2\sum_{k\in B_j}(2\pi
k)^4(f_k^2+\eps^2))^{1/2}}{4(2\eps\pi \kappa_{j+1})^2}.
$$
It is clear that
\begin{align*}
&Q_\eps\geq \frac{(8\sum_{k\in B_j}(2\pi k)^4)^{1/2}}{4(2\pi
\kappa_{j+1})^2}\geq
\frac{(8(2\pi)^4\kappa_{j+1}^4T_j/5)^{1/2}}{4(2\pi \kappa_{j+1})^2}\geq\sqrt{T_j/10}\\
&\sum_{k\in B_j}(2\pi k)^4(f_k^2+\eps^2)\le (2\pi \kappa_{j+1})^2(\|f'\|^2_{(j)}+\eps^2\sigma_j^2)
\le 3T_j^{-1}\sigma_j^2(\|f'\|^2_{(j)}+\eps^2\sigma_j^2)\\
&\qquad\qquad\qquad\qquad\quad\le \frac{3(1-h_{\kappa_j}^*)(\|f'\|^2_{(j)}+\eps^2\sigma_j^2)^2}{T_j\eps^2}.
\end{align*}
These inequalities combined with the identity $\eta_j=\sum_{k\in B_j}Y_j/(\|f'\|^2_{(j)}+\eps^2\sigma^2_j)$
yield the desired result.
\end{proof}

\begin{lemma}\label{lemh_k} Assume that $\varphi_j^2T_j\geq 24\log\eps^{-5}$ and $\varphi_j<1$.
Then
\begin{align*}
\Pb\bigg(\hat h_k^{JS}\in\bigg[h_k^*-\frac{2(1-h_k^*)\varphi_j}{(1-\varphi_j)}\,,h_k^*\bigg]\bigg)\ge 1-2\eps^5.
\end{align*}
\end{lemma}
\begin{proof}
Note that
\begin{equation}\label{h1}
\hat h_k^{JS}=\bigg(1-\frac{(1-h_k^*)(1+\varphi_j)}{1+\eta_j}\bigg)\1_{\{\eta_j>-h_k^*+\varphi_j(1-h_k^*)\}}.
\end{equation}
One checks that $\hat h_k^{JS}>h_k^*$ if and only if $\eta_j>\varphi_j$. Therefore, $\Pb(\hat h_k^{JS}>h_k^*)
=\Pb(\eta_j>\varphi_j)\leq \eps^4$.
Similarly,
\begin{align*}
\Pb\bigg(\hat h_k^{JS}<h_k^*-\frac{2(1-h_k^*)\varphi_j}{(1-\varphi_j)}\bigg)&=\Pb\Big(\varphi_j-(1+\varphi_j)h_k^*\le\eta_j\leq -\varphi_j\Big)\\
&\quad+\Pb(\eta_j<\varphi_j-(1+\varphi_j)h_k^*)\1_{\{h_k^*>\frac{2(1-h_k^*)\varphi_j}{(1-\varphi_j)}\}}\\
&\le \Pb\big(\eta_j\leq -\varphi_j\big),
\end{align*}
since $h_k^*>\frac{2(1-h_k^*)\varphi_j}{(1-\varphi_j)}$ if and only if $\varphi_j-(1+\varphi_j)h_k^*<-\varphi_j$. Therefore, using
Lemma~\ref{lemeta},
$$
\Pb\bigg(\hat h_k^{JS}\in\bigg[h_k^*-\frac{2(1-h_k^*)\varphi_j}{(1-\varphi_j)}\,,h_k^*\bigg]\bigg)\ge \Pb\big(|\eta_j|\le \varphi_j\big)\ge
1-2\eps^5
$$
and the assertion of the lemma follows.\end{proof}

\begin{lemma}\label{lem10}
For any positive $x$ verifying $ x^2\le T_\eps/5$ it holds
$$
\Pb\bigg(\bigg|\sum_{k=1}^{N_\eps}\tilde h_k(2\pi
k)^2(\xi_k^2-{\xi_k^*}^2)\bigg|\ge
\frac{12\sqrt{2}x}{\sqrt{T_\eps}}\;R^\eps[f,h^*]\bigg)\le
2Je^{-x^2}.
$$
\end{lemma}
\begin{proof}
Using (\ref{3Tj}), we get $(\sum_{k\in B_j}(2\pi k)^4\big)^{1/2}\le 3\sigma_j^2/\sqrt{T_j}$.
This inequality combined with the fact that  $\tilde h_k\leq h_k^*$, allows us to bound the
probability of the event of interest by
\begin{align*}
\Pb\bigg(\bigg|\sum_{k=1}^{N_\eps}&\tilde h_k(2\pi
k)^2(\xi_k^2-{\xi_k^*}^2)\bigg|\ge
\sqrt{32}\;x\sum_{j=1}^J h_{\kappa_j}^*\Big(\sum_{k\in B_j}(2\pi k)^4\Big)^\frac12\bigg)\\
&\le\sum_{j=1}^J \Pb\bigg(\tilde h_{\kappa_j}\bigg|\sum_{k\in B_j}
(2\pi k)^2(\xi_k^2-{\xi_k^*}^2)\bigg|\ge \sqrt{32}\,xh_{\kappa_j}^*\Big(\sum_{k\in B_j}(2\pi k)^4\Big)^{1/2}\bigg)\\
&\le \sum_{j=1}^J \Pb\bigg(\bigg|\sum_{k\in B_j} (2\pi
k)^2(\xi_k^2-{\xi_k^*}^2)\bigg|\ge \sqrt{32}\,x\Big(\sum_{k\in
B_j}(2\pi k)^4\Big)^{1/2}\bigg).
\end{align*}
The desired result follows now from Lemma~\ref{varsigma} and~(\ref{mnp}).
\end{proof}

\begin{lemma}\label{lem11}
For any $x>0$,
$$ \Pb\Bigg(\eps\bigg|\sum_{k=1}^{N_\eps} (\tilde h_k-h_k^*)(2\pi
k)^2\xi_kf_k\bigg|\ge 5
x\sqrt{\frac{\varphi_\eps}{T_\eps}}\;R^\eps[f,h^*]\Bigg) \leq
2Je^{-x^2}\!\!\!\!\!.
$$
\end{lemma}
\begin{proof}
Remark first that
$$ (\tilde h_{\kappa_j}-h_{\kappa_j}^*)^2\le
\hkapjsq\wedge 8(1-h_{\kappa_j}^*)^2\varphi_j^2 \le
\sqrt{8}{\varphi_jh_{\kappa_j}^*(1-h_{\kappa_j}^*)}
$$
for all $j=1,\ldots, J$. Set $Y_j=\eps\sum_{k\in B_j}(2\pi k)^2\xi_k
f_k$. The random variables $Y_1,\ldots,Y_J$ are independent zero
mean Gaussian with variance
$$
\esp[Y_j^2]=\eps^2\sum_{k\in B_j}(2\pi k)^4f_k^2\leq
3\eps^2\sigma_j^2T_j^{-1}\sum_{k\in B_j}(2\pi k)^2f_k^2.
$$
Therefore,
$\Pb(|Y_j|\geq \sqrt{6/T_j}\;x \eps\sigma_j\|f'\|_{(j)})\leq
2e^{-x^2}$
and consequently,
$$
\Pb\bigg(\sum_{j=1}^n|(\tilde h_{\kappa_j}-h^*_{\kappa_j})Y_j|\geq
\sqrt{\frac6{T_\eps}}\;x \eps\sum_{j=1}^J
\sqrt{4{\varphi_jh_{\kappa_j}^*(1-h_{\kappa_j}^*)}}\sigma_j\|f'\|_{(j)}\bigg)\leq
2Je^{-x^2}\!\!\!.
$$
To complete the proof, note that
\begin{align*}
\eps\sum_{j=1}^J
\sqrt{h_{\kappa_j}^*(1-h_{\kappa_j}^*)}\;\sigma_j\|f'\|_{(j)}&\leq
\frac12\sum_{j=1}^J (1-h_{\kappa_j}^*)\|f'\|_{(j)}^2+ \frac12\sum_{j=1}^J
h_{\kappa_j}^*(\eps \sigma_j)^2
\end{align*}
and the right side is bounded by $R^\eps[f,h^*]$.
\end{proof}

\begin{lemma}\label{lem12}
We have $\displaystyle\sum_{k=1}^{N_\eps} (\tilde h_k-h_k^*)(2\pi k)^2f_k^2\leq 4\varphi_\eps R^\eps[f,h^*]$.
\end{lemma}
\begin{proof}
The desired inequality is trivially fulfilled on $\mathcal E^c$, while on $\mathcal E$ we have
$0\le h_k^*-\tilde h_k\le 4\varphi_\eps(1-h_k^*)$, $\forall k\in B_j$, and hence
\hglue15pt$\displaystyle
\sum_{k=1}^\infty |\tilde h_k-h_k^*|(2\pi k)^2f_k^2\leq
4\varphi_\eps\sum_{k=1}^{N_\eps}(1-h_k^*)(2\pi k)^2f_k^2\leq 4\varphi_\eps R^\eps[f,h^*].
$
\end{proof}

\begin{lemma}\label{lemPhi'''}
Set $\hat X=\max_{1\le k\le
N_\eps}\sqrt{\xi_k^2(\theta)+{\xi_k^*}^2(\theta)}$. There exists an
event of probability at least $1-2\eps^4$ such that on this event,
for all $\tau\in\RR$, we have
\begin{align*}
\frac{|\Phi_\eps'''(\tau,\tilde h)|}{R^\eps[f,h^*]}\leq&\
\frac{12|\tau-\theta|\cdot\|f'\|^2}{\eps^2T_\eps}
+\frac{12\|f'\|\sqrt{\log\eps^{-8}}}{\eps T_\eps}
\\
&\ + \frac{16\sqrt{3}\;\pi N_\eps|\tau-\theta|\hat X\|f'\|}{\eps
\sqrt{T_\eps}}+4\pi N_\eps\hat X^2.
\end{align*}
\end{lemma}

\begin{proof}
Using (\ref{phi'1}), one checks that
\begin{align*}
\Phi_\eps'''(\tau,h)&=2\sum_{k=1}^\infty h_k(2\pi k)^3[(f_k+\eps\xi_k(\theta))^2-\eps^2\xi_k^*(\theta)^2]\sin[4\pi k(\tau-\theta)]\nonumber\\
&\quad-4\eps\sum_{k=1}^\infty h_k(2\pi k)^3
(f_k+\eps\xi_k(\theta))\xi_k^*(\theta)\cos[4\pi k(\tau-\theta)].
\end{align*}
Without loss of generality, we assume in the sequel that $\theta=0$. Then
\begin{align*}
\big|\Phi_\eps'''(\tau,\tilde h)\big|&\le 2\sum_{k=1}^{N_\eps} \tilde
h_k(2\pi k)^3 f_k^2|\sin(4\pi k\tau )|
+4\eps\sum_{k=1}^{N_\eps} \tilde h_k(2\pi k)^3|f_k\xi_k^*|\\
&\qquad+4\eps\sum_{k=1}^{N_\eps} \tilde h_k(2\pi k)^3\big|f_k\big[\xi_k\sin(4\pi k\tau)+\xi_k^*(\cos(4\pi k\tau)-1)\big]\big|\\
&\qquad+2\eps^2\sum_{k=1}^{N_\eps} \tilde h_k(2\pi k)^3\big|(\xi_k^2-{\xi_k^*}^2)\sin(4\pi k\tau)-2\xi_k\xi_k^*\cos(4\pi k\tau)\big|.
\end{align*}
Using the inequalities $|\sin(4\pi k\tau )|\le 4\pi k|\tau|$ and
$$
\sum_{k=1}^{N_\eps} \tilde h_k (2\pi k)^4 f_k^2\leq \sum_{j=1}^J h_{\kappa_j}^* \|f'\|_{(j)}^2 (2\pi \kappa_{j+1})^2\leq
\frac{3\|f'\|^2}{T_\eps} \sum_{j=1}^J h_{\kappa_j}^* \sigma_j^2,
$$
as well as the inequality $\eps^2\sum_{j=1}^J h_{\kappa_j}^* \sigma_j^2\le R^\eps[f,h^*]$, we get the desired bound for the first sum.
The bound on the second term is obtained using Lemma~\ref{lemsym}, the well known bound on the Laplace transform of a
Gaussian distribution and the inequality
\begin{align*}
\sum_{k=1}^{N_\eps} {h_k^*}^2(2\pi k)^6f_k^2&\leq \sum_{j=1}^J
\hkapjsq(2\pi \kappa_{j+1})^4\|f'\|_{(j)}^2\leq
\frac{9}{T_\eps^2}\sum_{j=1}^J\hkapjsq\sigma_j^4\|f'\|_{(j)}^2\\
&\leq \frac{9\|f'\|^2}{T_\eps^2}\max_j \hkapjsq\sigma_j^4\leq
\frac{9\|f'\|^2R^\eps[f,h^*]^2}{\eps^4T_\eps^2}.
\end{align*}
The bounds on the two remaining sums are obtained by combining the inequalities
\begin{align*}
|\xi_k\sin(4\pi k\tau)+\xi_k^*(\cos(4\pi k\tau)-1)|&\le |4\pi k\tau|\cdot|\xi_k\cos(4\pi k\tilde\tau)-\xi_k^*\sin(4\pi k\tilde\tau)|\\
&\leq 4\pi k|\tau| \hat X,\\
|(\xi_k^2-{\xi_k^*}^2)\sin(4\pi k\tau)-2\xi_k\xi_k^*&\cos(4\pi k\tau)|\leq \hat X^2
\end{align*}
with arguments similar to those used to bound the first two sums.
\end{proof}

\begin{lemma}\label{lem13} For any $x>0$, it holds
$$
\Pb\bigg(\Big|\Phi_\eps'(\theta,\tilde
h)-\Phi_\eps'(\theta,h^*)\Big|^2\ge 12\eps^2 x^2 \varphi_\eps
R^\eps[f,h^*] \bigg)\le 2e^{-x^2(1\wedge \sqrt{T_\eps/8x^2})}.
$$
\end{lemma}
\begin{proof} Let us denote $X_k=(2\pi k)(f_k+\eps\xi_k)\xi_k^*$.
According to (\ref{hatPhi'}),
\begin{align*}
\Phi_\eps'(\theta,\tilde h)-\Phi_\eps'(\theta,h^*)&
=\eps\sum_{j=1}^k(\tilde h_k-h_k^*)X_k.
\end{align*}
According to (\ref{e^tX}), for all $t\leq 1/(2\sqrt{2}\pi k\eps)$,
we have
$\esp[e^{tX_k}]\leq e^{t^2(2\pi k)^2(f_k^2+\eps^2)}$.
Thus the conditions of Lemma~\ref{lemsym} are fulfilled with
$\varrho_k=\tilde h_k-h_k^*$,
$\varrho_k^*=3\sqrt{h_k^*\varphi_j(1-h_k^*)}$, $T_k=1/(\sqrt8\,\pi
k\eps)$ and $g_k^2=(2\pi k)^2(f_k^2+\eps^2)$. Therefore,
$$
\Pb\bigg(\Big|\Phi_\eps'(\theta,\tilde
h)-\Phi_\eps'(\theta,h^*)\Big|^2\ge 12\eps^2 x^2\varphi_\eps
\sum_{k=1}^{N_\eps} (2\pi k)^2 h_k^*(1-h_k^*)(f_k^2+\eps^2)\bigg)\le
2e^{-x^2}
$$
for any $x>0$ verifying
$$
x^2\leq \frac{3\sum_{k=1}^{N_\eps} \varphi_{j(k)}(2\pi k)^2
h_k^*(1-h_k^*)(f_k^2+\eps^2)}{8\eps^2\max_j
\varphi_{j}h_{\kappa_j}^*(1-h_{\kappa_j}^*)(2\pi \kappa_{j+1})^2}.
$$
To complete the proof, it suffices to remark that
$$
\sum_{k=1}^{N_\eps} (2\pi k)^2 h_k^*(1-h_k^*)(f_k^2+\eps^2)= \sum_{j=1}^J
h_{\kappa_j}^*(1-h_{\kappa_j}^*)
[\|f'\|_{(j)}^2+\eps^2\sigma_j^2]\leq R^\eps[f,h^*]
$$
and
$$
\frac{3\sum_{k=1}^{N_\eps} \varphi_{j(k)}(2\pi k)^2
h_k^*(1-h_k^*)(f_k^2+\eps^2)}{8\eps^2\max_j
\varphi_{j}h_{\kappa_j}^*(1-h_{\kappa_j}^*)(2\pi
\kappa_{j+1})^2}\geq \frac{3\max_j
\varphi_{j}h_{\kappa_j}^*(1-h_{\kappa_j}^*)\sigma_j^2}{8\max_j
\varphi_{j}h_{\kappa_j}^*(1-h_{\kappa_j}^*)(2\pi \kappa_{j+1})^2}.
$$
Since $\sigma_j^2\geq T_j(2\pi \kappa_{j+1})^2/3$, the assertion of the lemma follows.
\end{proof}

\begin{lemma}\label{lemPhi'z} Let $\zeta=2\eps\|f'\|^{-2}\sum_{k} h_k^*(2\pi k)^2f_k\xi_k$.
For any event $\mathcal A$ verifying $\Pb(\mathcal A^c)=O(\eps^4)$, we have
$$
\esp[\Phi'_\eps(\theta,h^*)^2\zeta\1_{\mathcal
A}]=o(\eps^2R^\eps[f,h^*]).
$$
\end{lemma}
\begin{proof}
We have
\begin{align*}
\esp[\Phi'_\eps(\theta,h^*)^2\zeta]&=\eps^2\esp\bigg[\sum_{k=1}^{N_\eps}
(2\pi k)^2{h_k^*}^2(f_k+\eps\xi_k(\theta))^2\xi_k^*(\theta)^2\zeta\bigg]\\
&=2\eps^3\esp\bigg[\sum_{k=1}^{N_\eps} (2\pi k)^2 {h_k^*}^2f_k
\xi_k\zeta\bigg]=4\eps^4\|f'\|^{-2}\sum_{k=1}^{N_\eps} (2\pi k)^4 {h_k^*}^3 f_k^2\\
&\le 12T_\eps^{-1}\eps^4\sum_{j=1}^J {h_{\kappa_j}^*}^3
\sigma_j^2\leq 12T_\eps^{-1}\eps^2R^\eps[f,h^*].
\end{align*}
Using the Rosenthal inequality, one easily checks that
$\esp[\Phi'_\eps(\theta,h^*)^{2p}]=O(\eps^{2p})$ and
$\esp[\zeta^{2p}]=o({R^\eps}[f,h^*]^{p})$ for any integer $p>0$.
Therefore, the Cauchy-Schwarz inequality yields,
$$
\big|\esp[\Phi'_\eps(\theta,h^*)^2\zeta\1_{\mathcal A^c}]\big|\leq
o\big(\eps^2\sqrt{R^\eps[f,h^*]}\;\big)\sqrt{\Pb(\mathcal
A^c)}=o\big(\eps^4\sqrt{R^\eps[f,h^*]}\;\big)
$$
and the assertion of the lemma follows.
\end{proof}

\subsection{Lemma used in Theorem~\ref{Th2}}
We assume that $f\in\mathcal F_{\delta,\beta,L}(\bar f)$ with
$\bar f\in\mathcal F(\beta^*,L^*,\rho)$ and $\beta^*>\beta\geq\beta_*$. For the sake of completeness we give below
a suitable version of the Pinsker theorem \cite{PIN}.

\begin{lemma}\label{lemPinsk}
Set $\gamma_\eps=1/\log\eps^{-2}$, $W_\eps=\big(\frac{L}{\varepsilon ^{2}}\frac{(\beta+2)(2\beta+1)}
{(2\pi )^{2\beta}(\beta-1)}\big)^{1/(2\beta+1)}$ and define
$$
\lambda^*_{k}=\begin{cases}
1, & k\leq \gamma_{\varepsilon }W_{\varepsilon }, \\
\displaystyle{\bigg[ 1-\Big(\frac{k}{W_{\varepsilon }}\Big)
^{\beta-1}\bigg] _{+}}, & k>\gamma _{\varepsilon }W_{\varepsilon }.
\end{cases}
$$
The filter $\lambda_k^*$ satisfies
\begin{align*}
\sup_{f\in\mathcal F_{\delta,\beta,L}} R^\eps[f,\lambda^*]\le
(1+o(1))\,C(\beta,L)\,\eps^{\frac{4\beta-4}{2\beta+1}}.
\end{align*}
\end{lemma}
\begin{proof}
Set $v=f-\bar f$. Using the inequality $(\bar f_k+v_k)^2\le2z^{-1}\bar f_k^2+(1+z)v_k^2$, $\forall z\in[0,1]$, we obtain
\begin{align}
R^\eps[f,\lambda^*]&=\sum_{k>\gamma_\eps W_\eps}(2\pi k)^2(1-\lambda_k^*)^2(\bar f_k+v_k)^2
+\eps^2\sum_{k=1}^{\infty}(2\pi k)^2{\lambda_k^*}^2\nonumber\\
&\le 2z^{-1}\sum_{k>\gamma_\eps W_\eps}(2\pi k)^2(1-\lambda_k^*)^2\bar f_k^2+(1+z)R^\eps[v,\lambda^*].\nonumber
\end{align}
Since $\bar f\in \mathcal F(\beta^*,L^*,\rho)$, we have
\begin{align}\label{barf}
\sum_{k>\gamma_\eps W_\eps}(2\pi k)^2(1-\lambda_k^*)^2\bar f_k^2&\le
L^*(\gamma_\eps W_\eps)^{2-2\beta^*}=o(W_\eps^{2-2\beta}).
\end{align}
On the other hand, setting $\tilde\lambda^*_k=(1-(k/W_\eps)^{\beta-1})_+$,
$$
R^\eps[v,\lambda^*]-R^\eps[v,\tilde\lambda^*]\leq \eps^2\sum_{k\le \gamma_\eps W_\eps} (2\pi k)^2\leq 4\pi^2\eps^2(\gamma_\eps W_\eps)^3
=o(\eps^2W_\eps^3).
$$
Using the relation $W_\eps^{2-2\beta}=O(\eps^2W_\eps^3)$ and choosing $z=z_\eps$ appropriately, we get
$$
\sup_{f\in\mathcal F_{\delta,\beta,L}(\bar f)} R^\eps[f,\lambda^*]\le (1+o(1))\sup_{v\in\mathcal W(\beta,L)}
R^\eps[v,\tilde\lambda^*]+o(\eps^2W_\eps^3),
$$
where $\mathcal W(\beta,L)$ is the Sobolev ball $\{v:\sum_{k\ge 1}(2\pi k)^{2\beta}f_k^2\le L\}$. It then follows from
\cite[Thm.\ 1 and Example 1, p. 265]{BeLe} (with $\alpha=\beta-1$ and $\delta=3$) that
$\sup_{v\in\mathcal W(\beta,L)}R^\eps[v,\tilde\lambda^*]=C(\beta,L)\,\eps^{\frac{4\beta-4}{2\beta+1}}(1+o(1))$.
To conclude, it suffices to remark that
$o(\eps^2W_\eps^3)=o(\eps^{\frac{4\beta-4}{2\beta+1}})$.
\end{proof}

\subsection{Auxiliary general results}

\begin{lemma}\label{varsigma} Assume that $a_1,a_2,\ldots,a_n\in \RR$ and $\varsigma=\sum_{k=1}^na_k(\xi_k^2-{\xi_k^*}^2)$, where
$(\xi_1,\ldots,\xi_n,\xi_1^*,\ldots,\xi_n^*)$ is a zero mean
Gaussian vector with identity covariance matrix. For any $y\in
[0,\|a\|/\max_k a_k]$, it holds
$$
\Pb\bigg(\varsigma^2\ge 32y^2\sum_{k=1}^n a_k^2\bigg)\leq 2
e^{-y^2}.
$$
\end{lemma}
\begin{proof} Using the formula of the Laplace transform of a chi-squared distribution,
for any $t\in[-(\sqrt8 a_k)^{-1},(\sqrt8a_k)^{-1}]$, we get
$$
\esp[e^{a_kt(\xi_k^2-{\xi_k^*}^2)}]=\frac1{1-4a_k^2t^2}\leq
e^{8a_k^2t^2}.
$$
Applying \cite[Thm.\ 2.7]{PET} with $g_k=16a_k^2$ and
$x=\sqrt{32}y\|a\|$, we get the desired result.
\end{proof}

\begin{lemma}\label{lemsym}
Let $X_1,\ldots,X_n$ be independent symmetric random variables. Let
$\boldsymbol{\varrho}=(\varrho_1,\ldots,\varrho_n)$ be a random
vector satisfying $|\varrho_j|\leq \varrho_j^*$, $\forall
j=1,\ldots,n$ with some deterministic sequence
$(\varrho_j^*)_{j=1}^n$ and $\mathscr
L(\boldsymbol{\varrho}|X_j=x)=\mathscr
L(\boldsymbol{\varrho}|X_j=-x)$ for all $j\in\{1,\ldots,n\}$. If
$$
\esp[e^{tX_j}]\leq e^{t^2g_j^2}
$$
for some sequence
$(g_j)_{j=1,\ldots,n}$ and for $|t|\leq T_j$, then
$$
\Pb\bigg(\bigg|\sum_{j=1}^n\varrho_j X_j \bigg|^2\ge 4x^2{\sum_{j=1}^n{\varrho_j^*}^2g_j^2}\bigg)\leq 2e^{-x^2(1\wedge Q_nx^{-1})},\
\forall\, x>0
$$
where $Q_n=(\sum_{j=1}^n{\varrho_j^*}^2g_j^2)^{1/2}\min_j(T_j/\varrho_j^*)$.
\end{lemma}
\begin{proof} Set $Y_j=\varrho_j X_j$ and $\bar Y_j=\varrho_j^* X_j$.
For any $p_1,\ldots,p_n\in\NN^n$, the expectation $\esp[Y_1^{p_1}\cdot\ldots\cdot Y_n^{p_n}]$ vanishes
if at least one $p_j$ is odd. Therefore, $\esp[(\sum_j Y_j)^k]=0$ if $k$ is odd and
$\esp[(\sum_j Y_j)^k]\leq \esp[(\sum_j \bar Y_j)^k]$ if $k$ is even. Hence
\begin{align*}
\esp\Big[\exp\Big(t\sum_j Y_j\Big)\Big]&=\sum_{k=0}^\infty\frac{t^k\esp[(\sum_j Y_j)^k]}{k!}\le \sum_{k=0}^\infty\frac{t^{2k}
\esp[(\sum_j\bar Y_j)^{2k}]}{(2k)!}\\
&=\esp[e^{t\sum_j\varrho_j^*X_j}]\leq e^{t^2\sum_j{\varrho_j^*}^2g_j^2},\quad \forall \,|t|\leq \min_j
({T_j}/{\varrho_j^*}).
\end{align*}
According to the Markov inequality, for every $t>0$,
$$
\Pb\bigg(\bigg|\sum_{j=1}^n\varrho_j X_j \bigg|\ge y\bigg)\leq 2e^{-ty}\esp[e^{t\sum_jY_j}]\leq
2\exp\bigg(-ty+{t^2}\sum_{j=1}^n{\varrho_j^*}^2g_j^2\bigg).
$$
Setting $t=(x\wedge Q_n)/(\sum_j {\varrho_j^*}^2g_j^2)^{1/2}$
and $y^2=4x^2{\sum_{j}{\varrho_j^*}^2g_j^2}$ we get the desired result.
\end{proof}

\noindent\textbf{Aknowledgement.} We are thankful to the anonymous referee for the remarks that helped
to improve the presentation.


\end{document}